\documentclass[runningheads]{llncs}

\usepackage{epsfig}

\usepackage{amssymb,amsmath} 
\usepackage{latexsym}

\newcommand{\ga}{ \gamma}
\newcommand{\pp}{{\bf p}}
\newcommand{\qq}{{\bf q}}
\newcommand{\mm}{{\bf m}}
\newcommand{\nn}{{\bf n}}
\newcommand{\ee}{{\bf e}}
\newcommand{\ff}{{\bf f}}

\newcommand{\PP}{\mathcal{P}}
\newcommand{\TT}{{\bf T}}

\newcommand{\RR}{{\mathbb R}}

\begin{document}
\title{Pre-Symmetry Sets of 3D shapes}
\author{Andr\'{e} Diatta  and Peter Giblin}
\institute{Department of Mathematical Sciences, University of Liverpool,\\ Liverpool L69 7ZL, England,\\  {\tt adiatta@liv.ac.uk  pjgiblin@liv.ac.uk}}
\authorrunning{A. Diatta and P. Giblin}   
\date{}
\maketitle
        
\begin{abstract}
We show that the pre-symmetry set of a smooth surface in 3-space has the structure of the graph of a function from $\RR^2$ to $\RR^2$ in
many cases of interest, generalising known results for the pre-symmetry set of a curve in the plane. We explain how this function is obtained, 
and illustrate with examples both on and off the diagonal. There are other cases where the
pre-symmetry set is {\em singular}; we mention some of these cases but leave their investigation to another occasion.
\end{abstract}

\section{Introduction}\label{s:intro}

In this paper, we consider in some detail the structure of pre-symmetry sets in 2D and 3D. Recall that, given $M$, a
smooth closed curve in 2D or surface in 3D, the {\em pre-symmetry set} $\PP$ of $M$ is the closure of the set of pairs of distinct points
$(\pp, \qq)\in M\times M$ for which there exists a circle or sphere tangent to $M$ at \pp\ and at \qq. From the pre-symmetry
set it is not difficult to pass to the {\em symmetry set} which is the locus of centres of these circles or spheres,
together with the centres of the limiting circles. When $M$
is a plane curve, parametrized by the points of a circle $S^1$, we can consider $\PP$ as a subset of the torus $S^1\times S^1$,
represented in the plane by a square with opposite sides identified. Note that in any dimension $\PP$ is symmetric:
$(\pp,\qq)\in \PP$ if and only if $(\qq,\pp)\in\PP$; it follows that in 2D we can also regard $\PP$ as contained in $S^1\times S^1$
with symmetric pairs identified. This is called the {\em symmetric product} of two circles and is a M\"{o}bius band
in which the boundary of the band represents the `diagonal' points $\pp=\qq$. In the 3D case, with $M$ say topologically equivalent
to a 2-sphere, we could consider $\PP$ as a subset of $S^2\times S^2$, which is topologically a complex quadric surface,
or of the symmetric product which is topologically a complex projective plane\footnote{This is a classic result
and is proved by regarding $S^2$ as the Riemann sphere, that is complex numbers together with $\infty$,
and then associating with an unordered pair $(z_1,z_2)$ of elements of $S^2$ the unique quadratic polynomial
$az^2+bz+c$ with roots $z_1$ and $z_2$. The corresponding element of the complex projective plane is then
$(a:b:c)$ and the diagonal corresponds to the conic $b^2=4ac$. 
In fact the symmetric product of any closed 2-dimensional manifold with itself is known to be a
2-manifold, a surprising fact since one might imagine that the diagonal would cause singularities.}.

However we shall not be concerned here with global models of the pre-symmetry set, but rather with
local or multi-local models. Starting from a circle or sphere $S_0$ having $k\ge 2$ contact points with our curve or
surface $M$, the pre-symmetry set can be partitioned into strata: \\
$\bullet$ \ one stratum for each pair of distinct
contact points, $\pp_i, \pp_j$ say, chosen from the $k$ points.  Here, we consider circles or spheres
tangent to $M$ at points {\em at or close to} $\pp_i$ and $\pp_j$ (if $k>2$ then these circles or spheres
`lose contact' with $M$ close to the other $k-2$ contact points of $S_0$ and $M$). \\
$\bullet$ \ one stratum for each contact point $\pp_i$ between $S_0$ and $M$ which is of type
$A_3$ or higher---a vertex and its circle of curvature in 2D, or  a ridge point and the corresponding sphere 
 of curvature in 3D. This stratum of $\PP$ arises from circles or spheres which are tangent to $M$
at two points both of which are close to $\pp_i$.

Note that we are concerned here with the symmetry set and not the medial axis: the circles or
spheres do not have to remain inside $M$ (or more generally be maximal with respect to $M$).
We can consider the above strata separately in our investigations. For example, starting with
a sphere $S_0$ having contact $A_1$ at two points $\pp_1, \pp_2$ and $A_3$ 
at $\pp_3$ with a surface $M$, we will have four strata which correspond to contacts
close to $\pp_1,\pp_3; \ \pp_2, \pp_3; \ \pp_1, \pp_2;$ and $\pp_3, \pp_3$.

In what follows we shall concentrate on the strata:\\
$\bullet$ arising from {\em two} contact points: {\em off-diagonal} strata, or\\
$\bullet$ arising from {\em one} contact point: {\em on-diagonal} strata\\
since the other cases can be reduced as above to these two. 
We shall argue that the `correct' way to think of these strata of the pre-symmetry
set $\PP$ in 2D ($n=2$) or 3D ($n=3$), when they are non-singular,  
is as {\em the graph of a mapping from $\RR^n$ to $\RR^n$.}
In this way we can capture the structure in a uniform way, both for a single curve or surface and for
a 1-parameter family of such. The mapping in question will arise in a slightly different way 
for the on-diagonal and off-diagonal cases, but the principle is the same
for both. Several of the standard examples of mappings from the plane to the plane~\cite{solidshape,rieger}---fold,
cusp, lips, beaks, swallowtail---arise naturally in this context. See Figure~\ref{fig:projs}.

The 2D case is relatively well-known~\cite{kuijper-olsen}; in \S\ref{s:2D} we summarize some results and interpret
one of the cases as the graph of a mapping $\RR\to\RR$ in order to set the scene for the 3D case.
We also take the opportunity to prove a result which explains why only certain singularities
of the pre-symmetry set can occur on the diagonal. 

In \S\ref{s:3D} we turn to 3D.
The underlying calculations here are rather complicated and we suppress them in favour of 
giving details of the results. This is part of a larger investigation of all the symmetry sets and pre-symmetry 
sets in 3D, for surfaces and generic 1-parameter families. Figure~\ref{fig:projs} summarises the results from \S\ref{s:3D}, Figure~\ref{fig:fold} treats the example of fold maps.

\section{The pre-symmetry set in 2D} \label{s:2D}

Given a smooth plane curve $M$ with parametrization $\ga:S^1\to \RR^2$ the pre-symmetry set is contained
in the set of pairs defined by the equations
\begin{equation}\label{eq:2d-pre-SS}
g(s,t)=0,   \text{ where } g(s,t)=(\gamma(s)-\gamma (t))\cdot(\TT(s)\pm \TT(t)),
\end{equation}
where \TT\ stands for the unit tangent vector.
The zero set of $g$ also contains diagonal pairs $(s,s)$, and pairs where the
tangents are parallel, and these need to be excluded when finding the true pre-symmetry set.

\subsection{The pre-symmetry set at a diagonal point}\label{ss:diag2D}

Let us examine the equation (\ref{eq:2d-pre-SS}) close to a diagonal point, that is,
when the curve $M$ has a vertex. We take $M$ in local form having
a vertex at the origin:
\[ \ga(x)=(x,y) \ \mbox{where} \ y=f(x)=a_2x^2+a_4x^4+a_5x^5 + \ldots, \ a_2\ne 0,\]
with no $x^3$ term. Expanding (\ref{eq:2d-pre-SS}) with the $+$ sign we find that there are no solutions
close to the origin besides $s=t$, and for the $-$ sign the leading terms (of degree 5) factorise as
\[ 2a_2(a_2^3-a_4)(s+t)(s-t)^4.\]
For an {\em ordinary vertex} (contact with the tangent circle of type $A_3$ exactly), $a_2^3\ne a_4$ and the local structure
of the pre-symmetry set, apart from the diagonal, is a smooth transverse curve $s+t$ + h.o.t. = 0.

For $A_4$ contact we need to go to the next terms which are
\[ -a_2a_5(3s^2+4st+3t^2)(s-t)^4,\]
where $a_5\ne 0$. The quadratic form here is positive definite, so
apart from the diagonal term there is an {\em isolated point} at $s=t=0$. Note that there
{\em cannot} be a pair of real branches of the pre-symmetry set at such a diagonal point.
We can also see this by noting that at an $A_4$ point the curvature $\kappa$ does not have an extremum
(since $\kappa'=\kappa''=0, \kappa'''\ne 0$) and on an arc without an extremum of
curvature there can be no bitangent circles~\cite{yuille-leyton}.

The same analysis can be continued to $A_5$, where the leading term is
\[ 2a_2(2a_2^5-a_6)(s+t)(2s^2+st+2t^2)(s-t)^4,\]
having a single transverse branch $s+t$ + h.o.t. = 0. The factor $2a_2^5-a_6$ is zero
only for an $A_6$ singularity. We find: \\
{\em for $k$ odd, a circle having $A_k$ contact with the curve $M$ results in a single branch
of the pre-symmetry set transverse to the diagonal; and for $k$ even it results in an isolated point of the pre-symmetry set
on the diagonal.}

\subsection{Transitions visible on the pre-symmetry set}\label{ss:trans2D}

The transitions which occur on symmetry sets were classified in~\cite{bg86}; some of these are `visible' on the
pre-symmetry set. There has been considerable work on this~\cite{kuijper-olsen}. We only give examples
here, together with an alternative interpretation of one of the cases, in line with our work on 3D pre-symmetry sets
below.

\begin{figure}
\begin{center}
\vspace*{-0.3in}
 \psfig{file=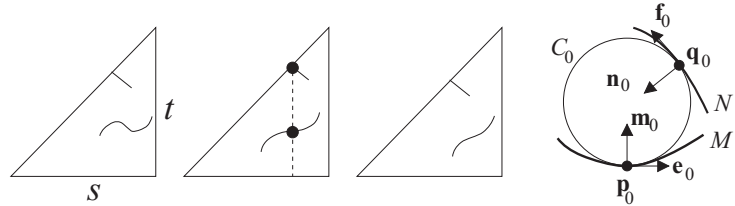,height=1.1in}
\end{center}
\vspace*{-0.3in}
\caption{Left three figures: an $A_1A_3$ transition seen in the pre-symmetry set; half of $S^1\times S^1$ is shown, the
$45^\circ$ line being the diagonal $s=t$. 
The transitional moment is in the centre,
with an inflexional tangent 
parallel to the $s$ parameter axis on one stratum of the pre-symmetry set.
There is also a diagonal point on the other stratum for the same value of $s$. Right: the
setup for analysing the $A_1A_3$ pre-symmetry set. The general bitangent circle under consideration
has contact points $\pp(s)$ close to $\pp_0$ and $\qq(t)$ close to $\qq_0$. Since $A_3$ occurs at the
`first' point $\pp_0$ and $A_1$ at the `second' point $\qq_0$ we could also write this as $A_3A_1$.}  
\label{fig:A1A3}
\end{figure}

\medskip\noindent
$\bullet$ \ $A_4$: the pre-symmetry set has an isolated point on the diagonal, growing to a closed loop
transverse to the diagonal in two places; compare~\S\ref{ss:diag2D}.\\
$\bullet$ \ $A_2^2$: two cases (1) moth, with the pre-symmetry set an isolated point growing into a closed loop,
and (2) nib, with the pre-symmetry set a transverse crossing of two branches separating two ways. Note that {\em off}
the diagonal, there is the possibility of both isolated point and transverse crossing on the pre-symmetry set,
in contrast with the situation {\em on} the diagonal described in~\S\ref{ss:diag2D}.\\
$\bullet$ \ $A_1A_3$: The pre-symmetry set has two strata, one on-diagonal and one off-diagonal. The off-diagonal
stratum has an inflexion parallel to one parameter axis and the line through the inflexion parallel to the
other parameter axis passes through the diagonal point. See Figure~\ref{fig:A1A3}.

\vspace*{-0.2in}

We consider in more detail the case of $A_1A_3$, as a preparation for the 3D case. 
Let $M, N$ be smooth arcs with arclength parameters $s, t$,  general points
$\pp=\pp(s), \qq=\qq(s)$ and curvatures $\kappa(s), \lambda(s)$ respectively.
Suppose that a circle $C_0$ is tangent to $M$ with contact $A_3$ at the point $\pp_0=\pp(0)$,
 and tangent to $N$ with contact $A_1$ at the point $\qq_0=\qq(0)$.
Write \ee, \mm\ for the unit tangent and normal at \pp\  and \ff, \nn\ for the unit tangent and normal at the point \qq.  Finally write
$r$ for the radius of a bitangent circle with contact points
\pp\ and \qq.  See Figure~\ref{fig:A1A3}.
 Since the circle $C_0$ has  $A_3$ contact with
$M$ at $s=0$, this point is a vertex of $M$ and we have $r\kappa=1$ and $d\kappa/ds=0$ at $s=0$.
 
We take as the defining
equation of the pre-symmetry set $F(s,t)={\bf 0}$ where
\begin{equation}
F(s,t)= \pp + r\mm - \qq - r\nn.
\label{eq:vector-preSS-2D}
\end{equation}
The Jacobian matrix of the mapping $F:\RR^3\to\RR^2$ is the $3\times 2$ matrix with columns the vectors
$\ee(1-r\kappa), \ \ \ -\ff(1-r\lambda), \ \ \ \mm-\nn$. Now,  using suffix 0 to mean
evaluation at $s=0$ or $t=0$, $\kappa_0=1/r, \ \lambda_0\ne 1/r$ so at $s=t=0$ the first of the column
vectors is zero and the others
are nonzero  multiples of $\ff_0, \ \mm_0-\nn_0$. The first is tangent to the circle $C_0$ and the second is parallel to
$\pp_0-\qq_0$ which is a chord of the circle; hence the two vectors are independent and we deduce from the implicit function
theorem that (\ref{eq:vector-preSS-2D}) has solution with $t$ and $r$ smooth functions of $s$ near $s=0$.

Now consider $t$ and $s$ as these functions of $s$ and write $t'$ for $dt/ds$, etc. Regarding (\ref{eq:vector-preSS-2D}) as
an {\em identity} in $s$ we can differentiate with respect to $s$ and get
\[ \ee(1-r\kappa)+r'\mm - \ff(1-r\lambda)t' - r'\nn = 0.\]
At $s=0$ the first term is zero and again $\mm-\nn, \ff$ are independent so we deduce $r'_0=t'_0=0$.
Differentiating again with respect to $s$, putting $s=0$ and keeping only the terms which are nonzero at $s=0$ gives
$r''(\mm-\nn)=\ff(1-r\lambda)t''$, so that $r''_0=t''_0=0$. Thus both $r$ and $t$, as functions of $s$, have degenerate critical
points. In fact differentiating again shows that for exactly $A_3$ at $s=0$ and $A_1$ at $t=0$ we have the next
derivatives of $r$ and $t$ nonzero. The fact that $t$ has $t'_0=t''_0=0$ is illustrated in Figure~\ref{fig:A1A3} where
the stratum of the pre-symmetry which is off-diagonal, namely the graph of $t$, has an {\em inflexion}.  Furthermore
with a little more trouble we can show that in a generic family of curves this {\em $A_2$ singularity of the function $t$}
is versally unfolded, that is it behaves as shown in Figure~\ref{fig:A1A3} with two critical points vanishing in the inflexion.

\begin{figure}
\begin{center}
\vspace*{-0.3in}
 \psfig{file=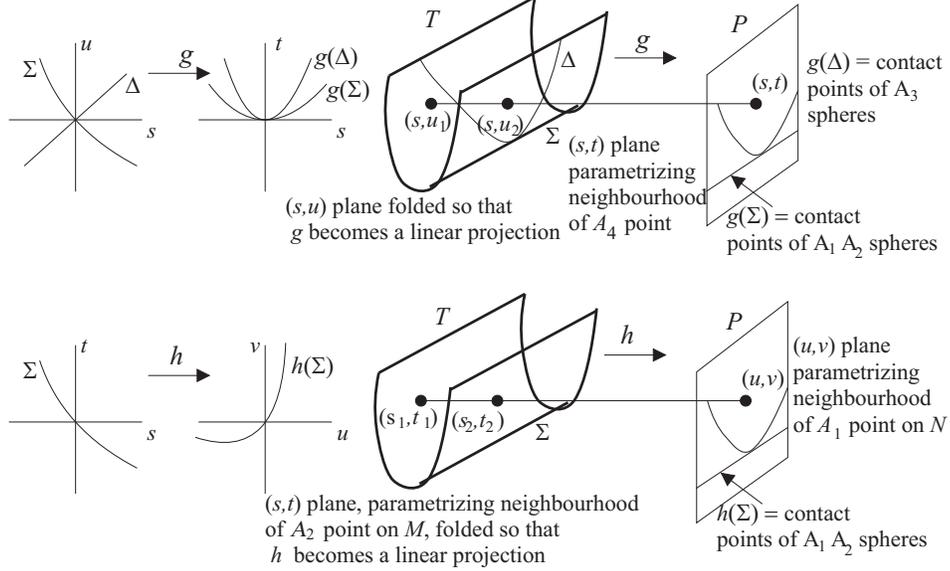,height=3in}
\end{center}
\vspace*{-0.3in}
\caption{Schematic diagrams of two of the mappings which we consider, both equivalent to folds. Above, Case $A_4$: 
$g$ maps a parameter plane $T$ to a neighbourhood $P$ of
the $A_4$ point. For a point $(s,t)$ `above' the fold line $g(\Sigma)$ on $P$, there are two bitangent spheres, with contact
at $(s,u_1), (u_1, v(s,u_1)$ and at $(s,u_2), (u_2, v(s, u_2))$. Below, Case $A_2A_1$: for a point $(u,v)$ `above' the
 fold line $h(\Sigma)$ on $P=$ neighbourhood of the $A_1$ point, there are bitangent spheres with contact at
$(s_1,t_1), (u,v)$ and at $(s_2,t_2), (u,v)$.}  
\label{fig:fold}
\end{figure}

\begin{figure}
\vspace*{-0.3in}
\begin{center}
 \psfig{file=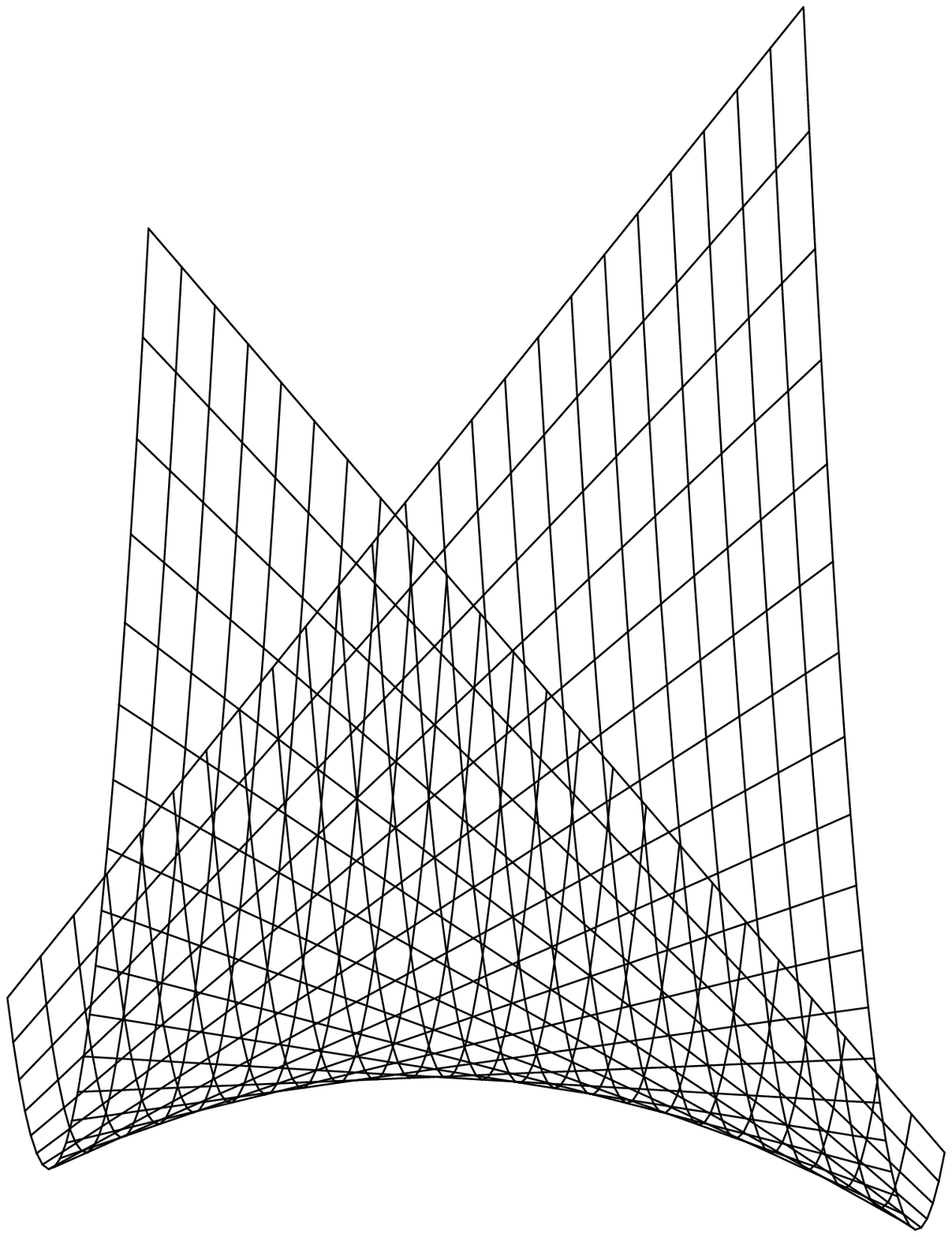,height=1.2in}
 \psfig{file=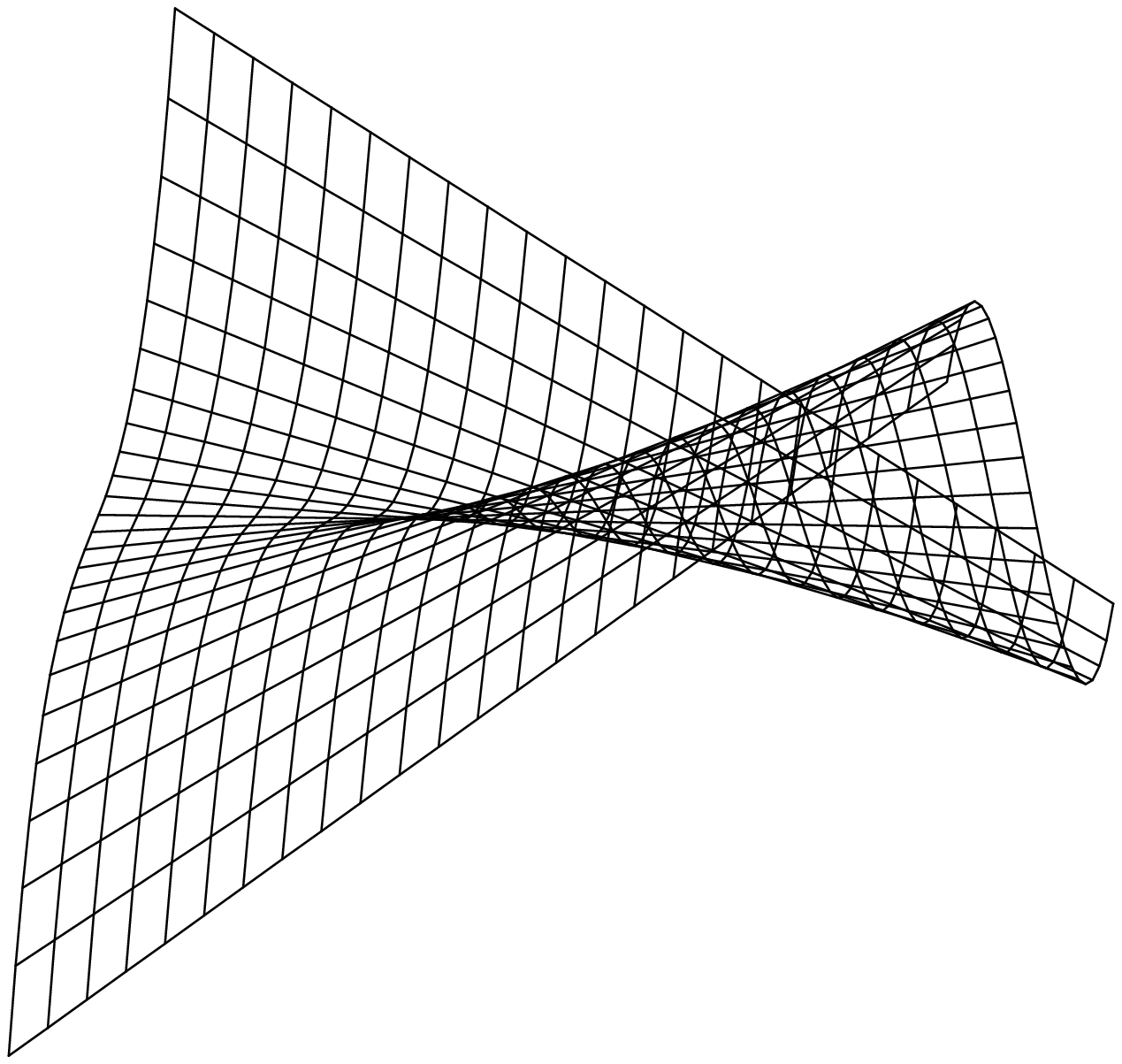,height=1.2in}
 \psfig{file=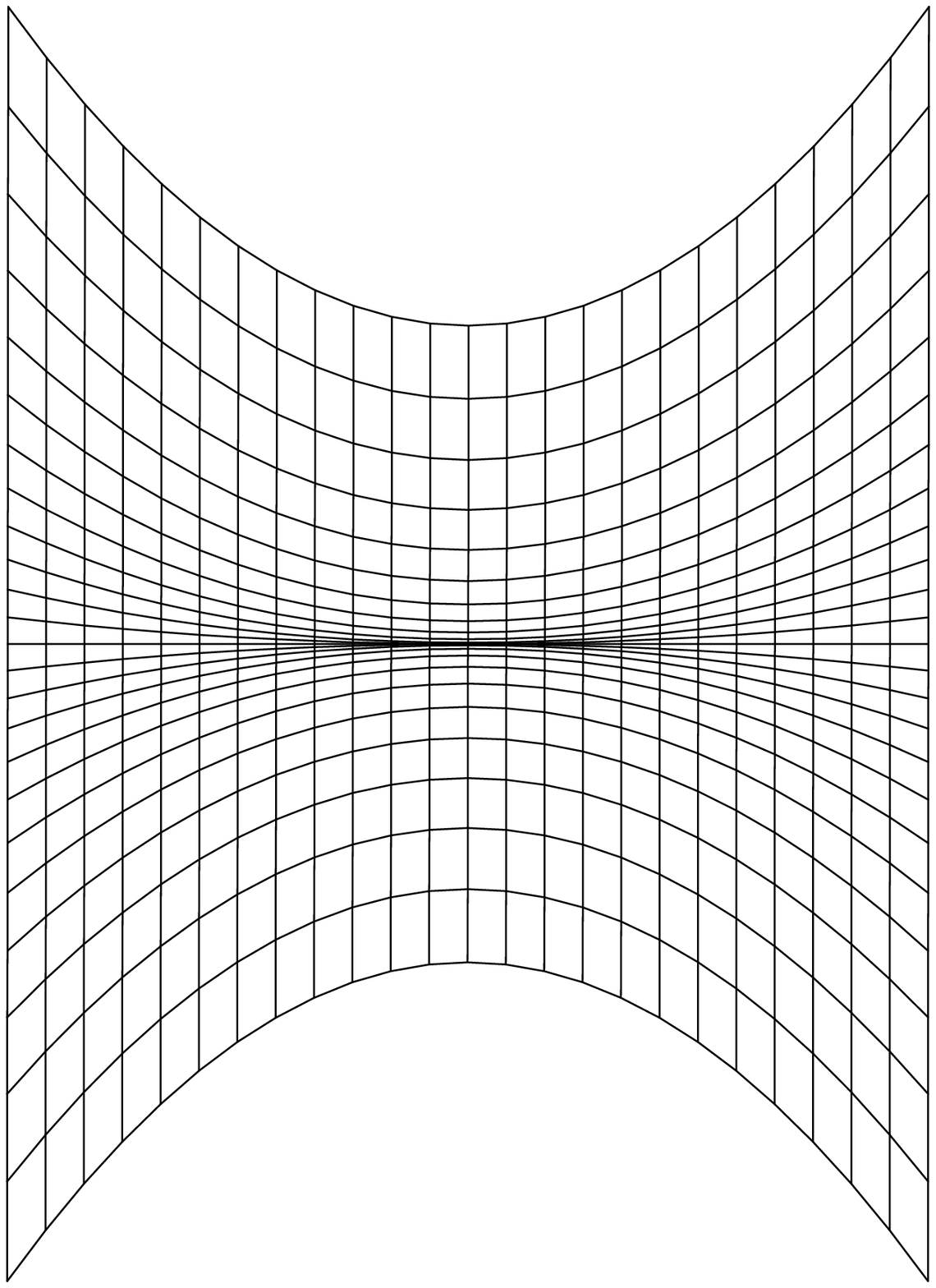,height=1.2in}
 \psfig{file=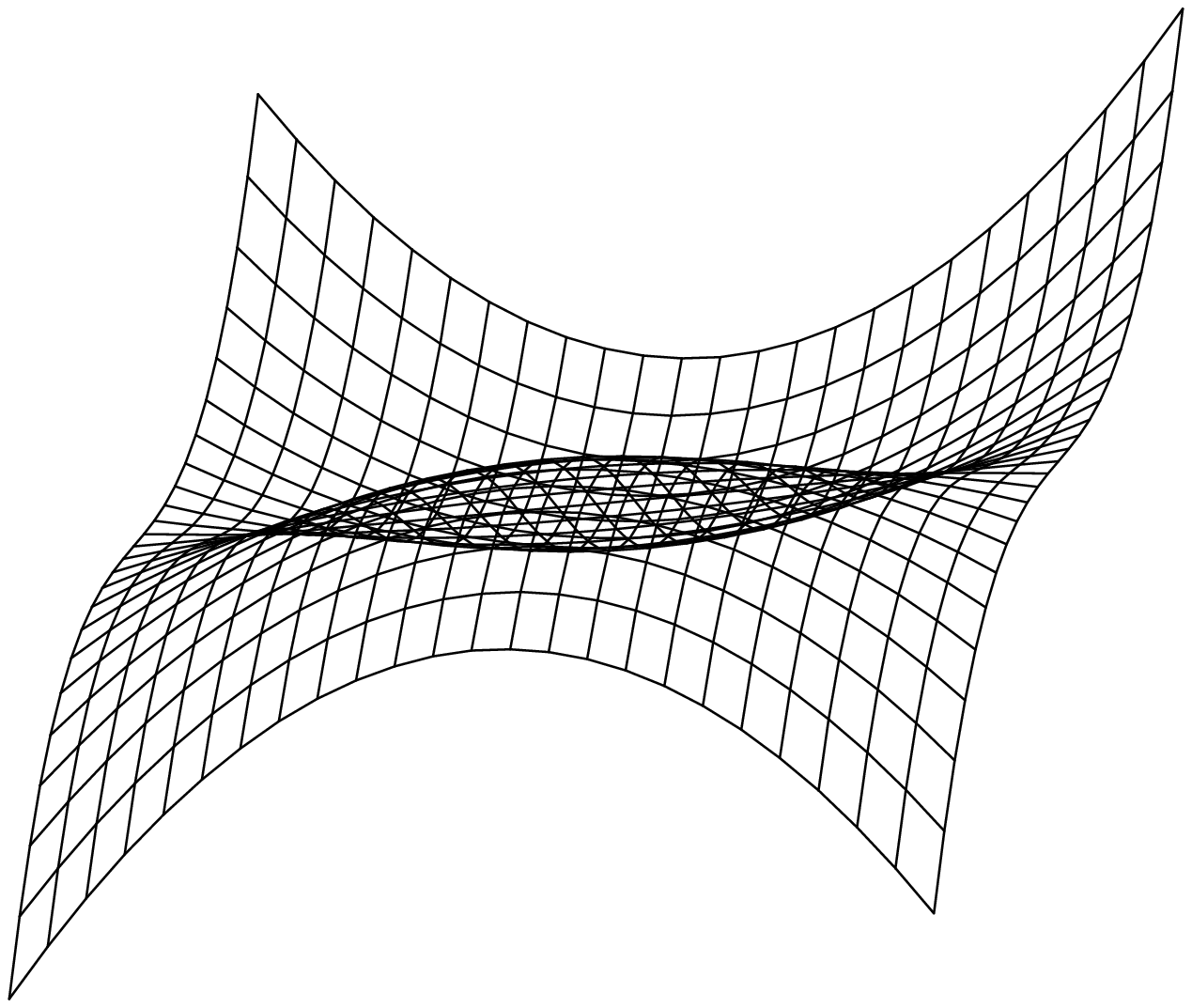,height=1.2in}
\end{center}
\vspace*{-0.2in}
\hspace*{0.5in} (1)\hspace*{1in} (2) \hspace*{0.7in}(3) \hspace*{1in}(3$'$)
\begin{center}
 \psfig{file=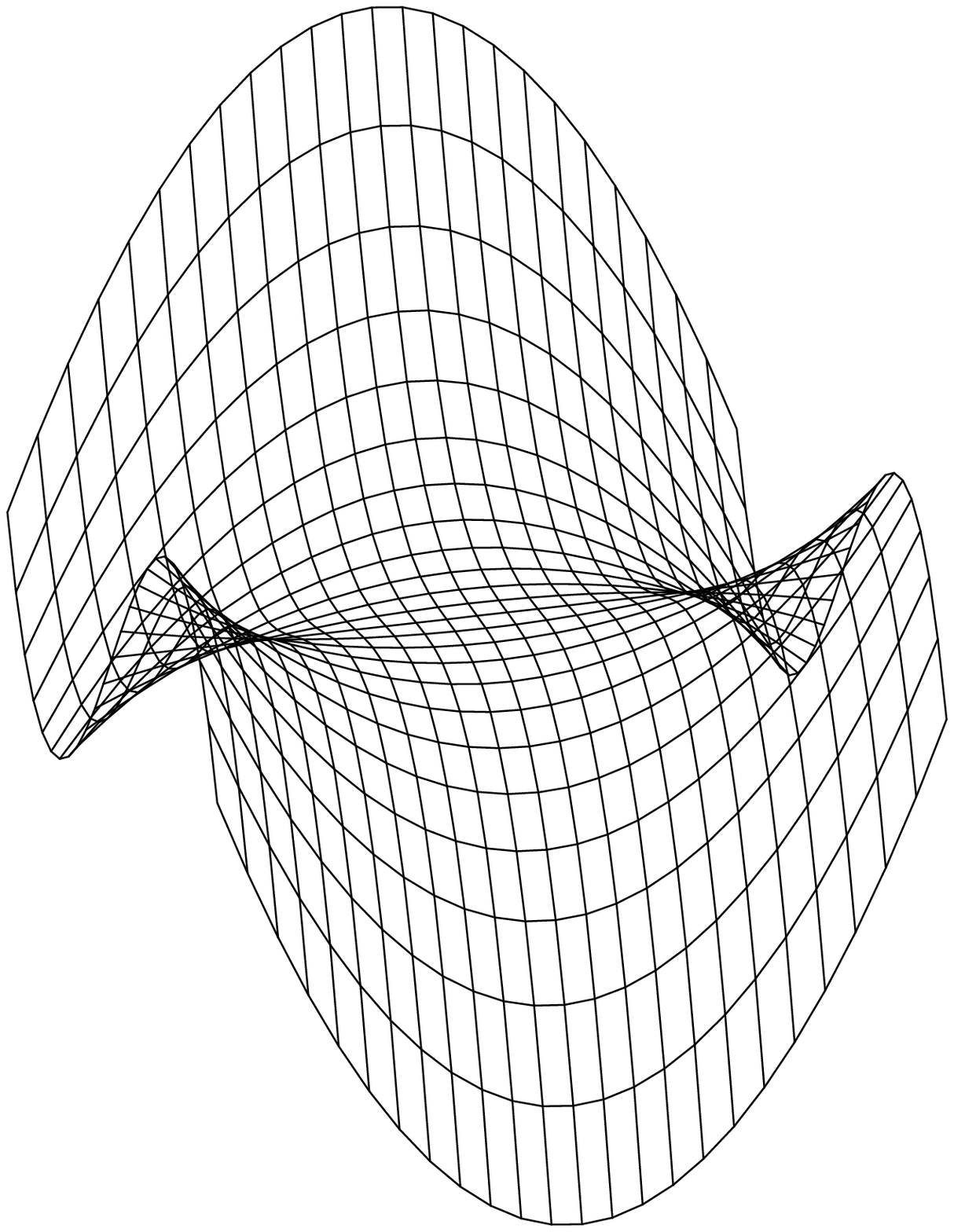,height=1.2in}
 \psfig{file=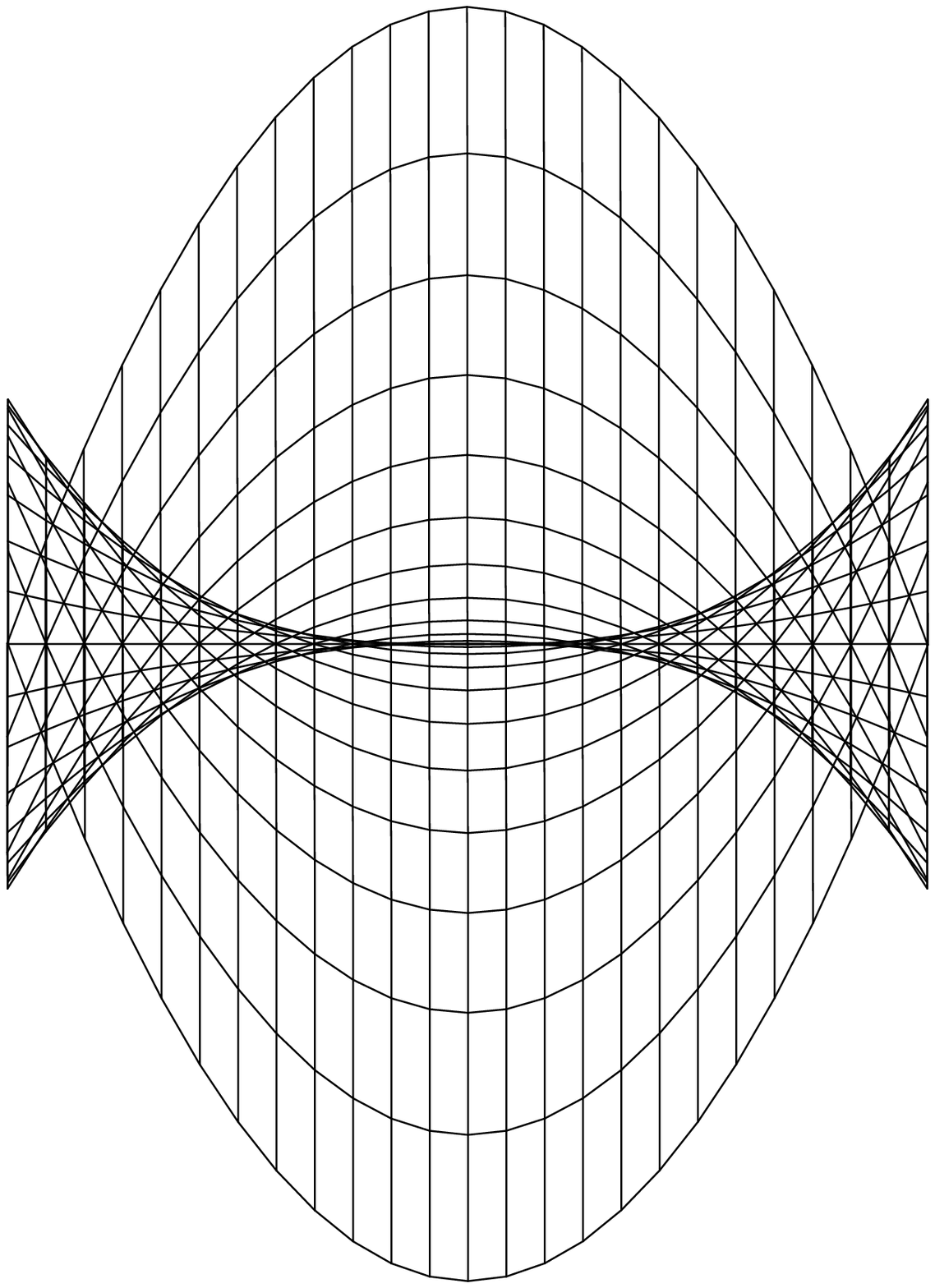,height=1.2in}
 \psfig{file=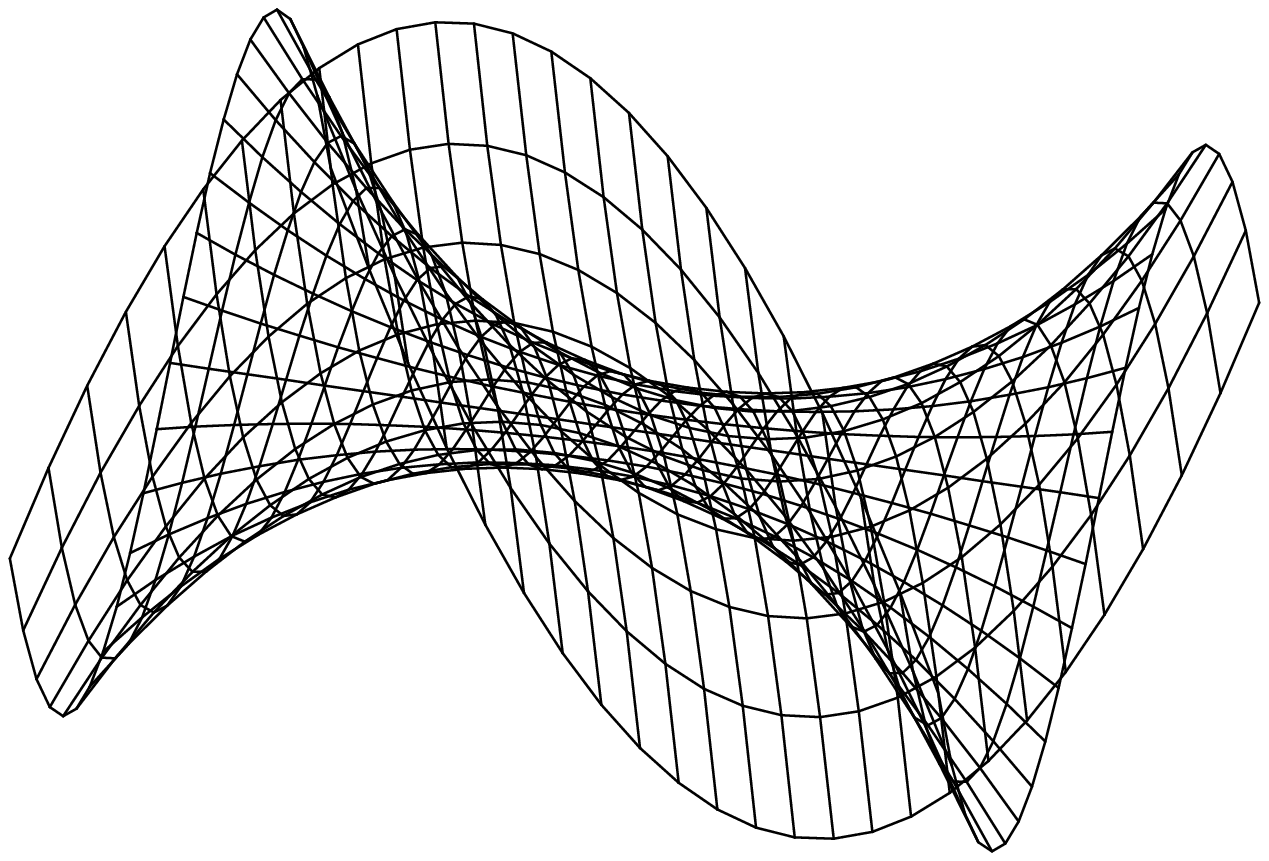,height=0.8in}
 \psfig{file=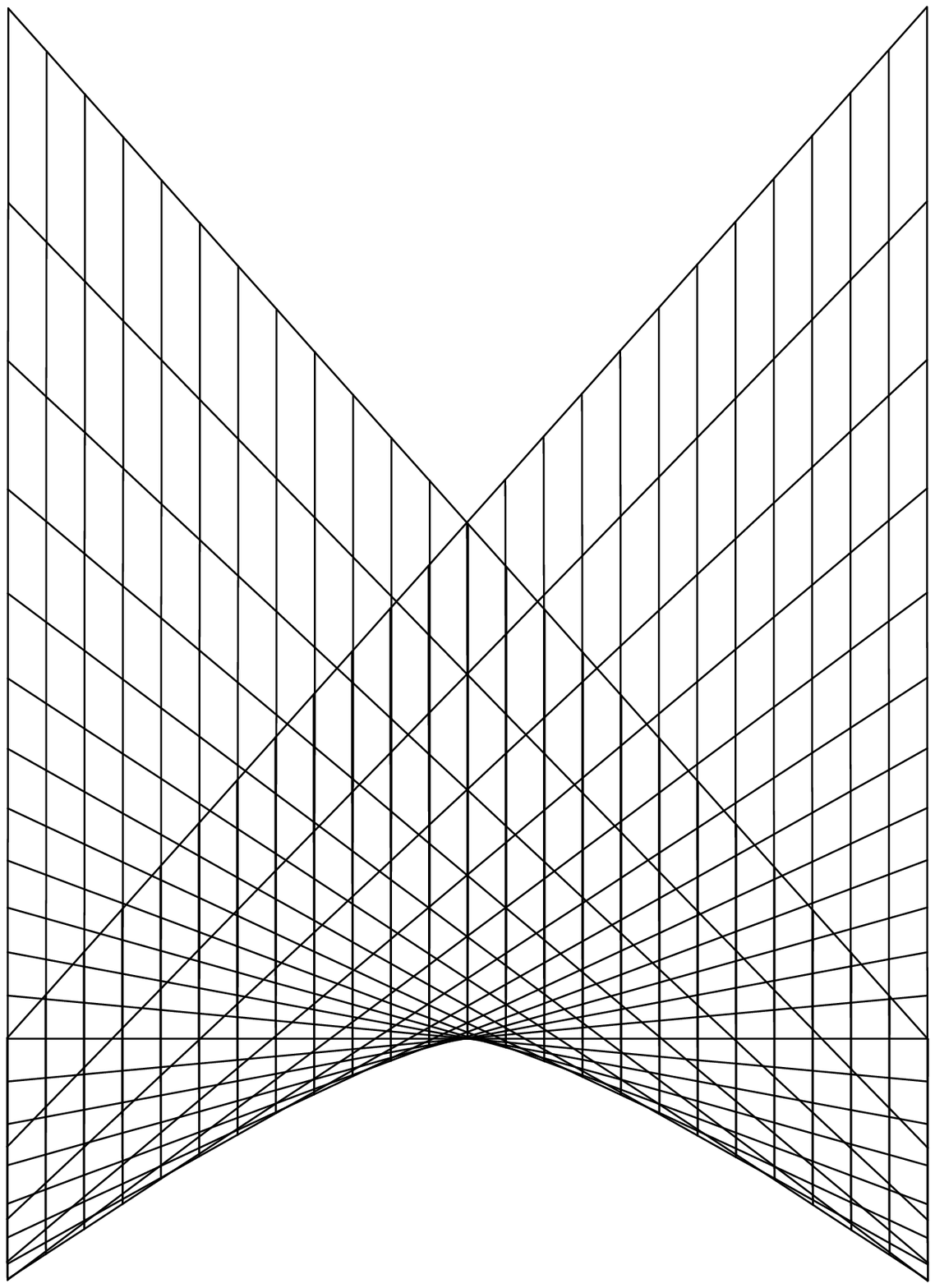,height=0.9in}
 \psfig{file=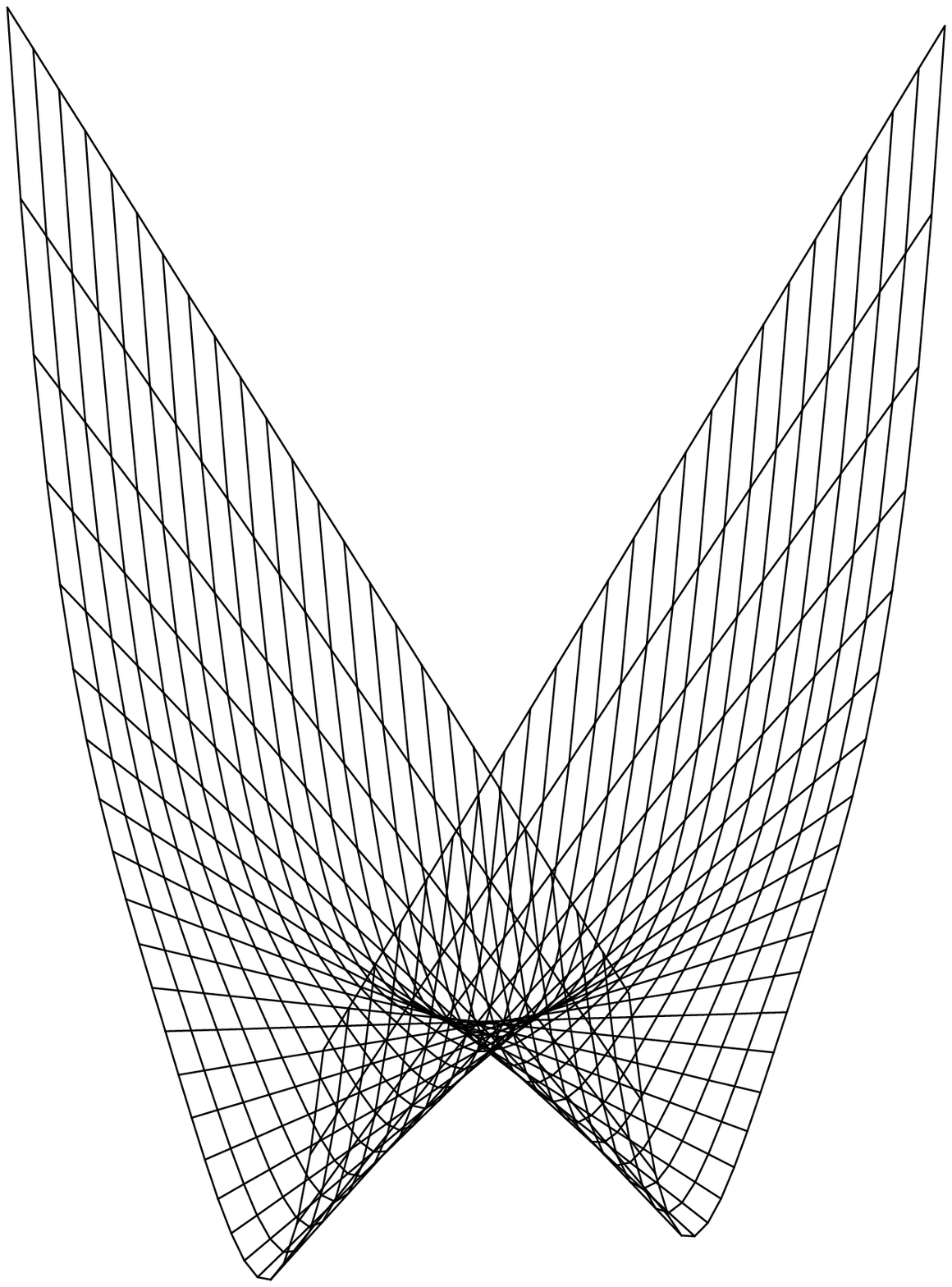,height=1.2in}
\end{center}
\vspace*{-0.1in}
\hspace*{0.3in} (4)\hspace*{0.8in} (4$'$) \hspace*{0.8in} (4$''$)\hspace{0.8in} (5) \hspace*{0.5in} (5$'$) 
\caption{Standard projections of a curved surface, say $T$ in $\RR^3$  to the plane $P$ of the paper. (1) fold; (2) cusp; (3) lips with unfolding $(3')$; $(4')$ beaks with
unfolding (4) and $(4'')$; and (5) swallowtail with unfolding $(5')$. These have the following
interpretations in the context of this paper. (1) Fold (compare Figure~\ref{fig:fold}): (i) Prop.~\ref{prop:AkA1cases}(b),  $T$ represents a neighbourhood
of an $A_2$ point on $M$ and $P$ represents a neighbourhood of an $A_1$ point on $N$ for an $A_2A_1$
singularity. Points of $P$ `below' the fold line have no corresponding  bitangent spheres on $T$ while
points of $P$ `above' the fold line have two. (ii) Prop.~\ref{prop:ondiagonal}(b), $A_4$. Now, $T$ represents the parameter plane of $s$ and $u$
while $P$ represents the surface $M$. Thus points of $M$ `below' the fold line are not contact points
for a bitangent sphere while points `above' are contact points for two such spheres. 
(2) Cusp:  Prop.~\ref{prop:AkA1cases}(c), generic $A_3A_1$ (fin point) only. The two fold lines ending in the cusp give points of $M$ (represented
by the curved surface $T$ of the figure) which are contact point of $A_2A_1$ spheres, as in (1). (3) Lips: (i) Prop.~\ref{prop:AkA1cases}(d): one
of the $A_3A_1$ transitional cases;  (ii) Prop.~\ref{prop:ondiagonal}(c), $A_5$.
 $(3')$ shows the situation immediately `after' the transition, with two cusps, representing generic $A_3A_1$ points
as in (2). (4) Beaks: the other $A_3A_1$ transition, Prop.~\ref{prop:AkA1cases}(d), only. This time two $A_2A_1$ curves
`survive' the transition though the $A_3A_1$ points (cusps as in (2)) are annihilated. (5) Swallowtail:  $A_4A_1$ transition, Prop.~\ref{prop:AkA1cases}(e), only.}
\label{fig:projs}
\end{figure}

We conclude:

\begin{proposition} At an $A_1A_3$ singularity the off-diagonal stratum of the pre-symmetry set has the structure of the
graph of a function with a  critical point of type $A_2$, that is a critical point equivalent to $t=s^3$, an ordinary
inflexion. In a generic family of curves the pre-symmetry set has  two ordinary critical points `before' the
$A_1A_3$ transition and none `after'.
\label{prop:A3A1-2D}
\end{proposition} 
In the 3D situation we shall also identify the pre-symmetry set as the graph of a function of a well-defined type. See illustrations in Figure~\ref{fig:fold} and Figure~\ref{fig:projs}.

\section{The pre-symmetry set in 3D}\label{s:3D}

In this section we shall investigate the strata of the pre-symmetry set of a smooth closed surface as in \S\ref{s:intro}.
We shall split into the off-diagonal and the on-diagonal cases; the results are strikingly similar but the details
are different. In fact we shall suppress most of the underlying mathematical calculations, which are similar
to those in \S\ref{s:2D} but more complicated since we are dealing here with surfaces instead of curves.

There are a great many cases of transitions on the symmetry set in 3D. Bogaevsky~\cite{bogaevsky} has determined
the transitions on the 3D medial axis (see also \cite{giblin-kimia-pollitt-3Dtransitions}) and in unpublished work has given
a list of all possible transitions on the 3D symmetry set, some of which have been investigated by Pollitt~\cite{pollitt-thesis}.
In this paper we shall determine the pre-symmetry set for a small number of representative cases; more mathematical detail
and a larger number of cases will appear elsewhere.

For the symmetry set of a generic surface $M$ in $\RR^3$ we consider here the following singularities:\\
\centerline{$A_3, \ A_4, \ A_1A_2, \ A_1A_3$,} \\
and for symmetry sets occurring in generic 1-parameter families of surfaces  we consider\\
\centerline{$A_1A_3$ transitions, $A_1A_4, \ A_5$.}
These results can also be used to determine other cases with three contact points,
such as $A_1^2A_2$ and $A_1^2A_3$. The cases $A_2^2, A_2A_3$ can be treated similarly.

\vspace*{-0.5in}

\subsection{The off-diagonal case}\label{ss:offdiag3D}
The basic setup is two pieces of smooth surface and a given sphere tangent to both of them at
known points. We then seek nearby points such that there is a sphere tangent at these nearby points.
To put the matter more precisely, consider two pieces of surface $M$ and $N$ and chosen points 
$\pp_0, \qq_0$ such that there is a sphere $S_0$ tangent to $M$ at $\pp_0$ and to $N$ at $\qq_0$. The order of tangency
at either point might be $A_k$ for $k$ an integer such that the total tangency is consistent with either a single
generic surface or a surface in a generic 1-parameter family of surfaces. For example $A_3$ at $\pp_0$ and
$A_1$ at $\qq_0$ would be a typical case. (The case of $D_4$ contact at an umbilic $\pp_0$ and $A_1$ contact elsewhere
is also generic in a 1-parameter family, and can be treated similarly.)

We take local parameters $(s,t)$ on $M$ near $\pp_0$  and $(u,v)$ on
$N$ near $\qq_0$, where $s=t=0$ at $\pp_0$ and $u=v=0$ at $\qq_0$, general points
on $M$ and $N$ being denoted \pp\ and \qq. 
Then the parameter space of points close to the pair $(\pp_0, \qq_0)$ is $\RR^4$, with coordinates
$(s,t,u,v)$. We seek the set of such points in $\RR^4$ such that there is a sphere tangent to $M$ at the point with
parameters $(s,t)$ and tangent to $N$ at the point with parameters $(u,v)$. This set is the local pre-symmetry set 
for the pair of surfaces $M, N$. As in \S\ref{ss:trans2D} we shall in many cases identify the pre-symmetry set with
the graph of a mapping of known type.

Let $r$ denote the radius function: the radius of the
given sphere $S_0$ is $r_0$.  Further,
let \mm\ be the unit normal at \pp, oriented so that for $\pp_0$ the normal passes through the centre of the given sphere,
and let \nn\ be the unit normal at \qq, oriented similarly.
Then the condition for a pre-symmetry set point is 
\begin{equation}
 G(s,t,u,v,r)=0 \text{ where }
G(s,t,u,v,r) = \pp + r\mm - \qq - r\nn.
\label{eq:pre-SS}
\end{equation}
In this equation, which is three scalar equations, there are five `unknowns' $s,t,u,v,r$ 
so we expect a two-dimensional solution. Note in particular that, on $G=0$,
the vectors $\pp-\qq$ and $\mm-\nn$ are parallel.

For the purpose of calculation we shall assume that
neither $\pp_0$ nor $\qq_0$ is an umbilic point of the corresponding surface. This enables us to use
principal directions for our coordinates around these points. Note that umbilics are isolated so we are losing only a finite
number of bitangent spheres.  Thus we shall assume the $s=$ constant and $t=$ constant curves on $M$ are, near
$\pp_0$, principal curves, and that for $s=0$ and for $t=0$ these curves are unit speed; and similarly for $N$.
Let $\ee_1, \ee_2$ be principal directions at \pp\ and $\ff_1, \ff_2$ principal directions at \qq.   Finally we
use $\kappa_1, \kappa_2$ for the principal curvatures of $M$ and $\lambda_1, \lambda_2$ for the principal
curvatures of  $N$.

The Jacobian matrix of $G$ at $(0,0,0,0,r_0)$ 
has the following vectors for its columns, using the fact that differentiating  the normal
in a principal direction \ee\ produces the corresponding principal curvature times \ee:
\begin{equation}
 (1 - r\kappa_1)\ee_1, \ \ 
 (1 - r\kappa_2)\ee_2, \ \ 
 -(1 - r\lambda_1)\ff_1, \ \ 
 -(1 - r\lambda_2)\ff_2, \ \
 \mm - \nn 
\label{eq:columns}
\end{equation}
Note that  $\mm-\nn$ is parallel to the chord joining $\pp$ and $\qq$, by (\ref{eq:pre-SS}).

The implicit function theorem gives us:
\begin{proposition}
Suppose that the contact at $\qq_0$ is of type $A_1$, which is the same as saying that
$r\ne 1/\lambda_1$ and $r\ne 1/\lambda_2$. Then $u,v,r$ are smooth functions of
$s$ and $t$ on the set $G=0$ and in particular the pre-symmetry set \\
$\{ (s,t,u,v): G(s,t,u,v,r)=0 \ \mbox{for some} \ r\}$\\
is the graph
of a smooth function $g:(s,t)\to (u,v)$.
\label{prop:AkA1}
\end{proposition}

The following result gives the nature of the mapping $g$ for several cases. The role of $g$ is that,
for any point \pp\ near $\pp_0$, $g$ gives a point \qq\ near $\qq_0$ such that there is a bitangent
sphere tangent to $M$ at \pp\ and to $N$ at \qq. Thus {\em every} point near $\pp_0$ is a possible
first point of contact but only points of the image of $g$ are possible second points of contact.
Furthermore, a given \qq\ may be the image of several points \pp\ under $g$. We illustrate these
properties in Figure~\ref{fig:projs} and Figure~\ref{fig:fold}.

We explain the
terms used in the Proposition after the statement.
\begin{proposition} 
\begin{enumerate}
\item[{\rm (a)}] In the $A_1A_1$ case, $g$ is a local diffeomorphism;
\item[{\rm (b)}] in the $A_2A_1$ case\footnote{We write $A_2A_1$ rather than $A_1A_2$ simply because the
`first' point $\pp_0$ is an $A_2$ and the `second' point $\qq_0$ is an $A_1$.}, $ g$ is a fold mapping;
\item[{\rm (c)}] in the {\em generic} $A_3A_1$ case (`fin point' on the medial axis), $g$ is a cusp mapping;
\item[{\rm (d)}] in the {\em transitional} $A_3A_1$ cases, $g$ is a lips or beaks mapping;
\item[{\rm (e)}] in the $A_4A_1$ case, $g$ is a swallowtail mapping.
\end{enumerate}
\label{prop:AkA1cases}
\end{proposition}
Part (c) here is to be compared with Proposition~\ref{prop:A3A1-2D} where the inflexion on the pre-symmetry set
in the  $A_1A_3$  case was interpreted as the graph of a mapping $\RR\to\RR$ of type $A_2$. Here we interpret
the structure of the pre-symmetry set as the graph of a cusp mapping $\RR^2\to\RR^2$. Note the distinction between
a {\em generic} $A_3A_1$, which occurs on a single smooth surface and a {\em transitional} $A_3A_1$ which
is one of the two forms identified by Bogaevsky~\cite{bogaevsky} in his classification of transitions on the 3D medial axis.

In the above Proposition, we use the standard names of classes of maps $\RR^2\to\RR^2$. These maps have
the following `normal forms'~\cite{rieger} up to smooth changes of coordinates in the two copies of $\RR^2$:\\
(a) local diffeomorphism: $(x,y)\to (x,y)$\\
(b) fold mapping: $(x,y)\to (x,y^2)$\\
(c) cusp mapping: $(x,y)\to (x,xy+y^3)$\\
(d) lips or beaks mapping: $(x,y)\to (x,x^2y\pm y^3)$ (+ lips, $-$ beaks)\\
(e) swallowtail mapping: $(x,y)\to(x,xy+y^4)$.

The different cases are separated by examining the critical set $\Sigma$ of $g$, that is the set of points
$(s,t)$ for which (using suffices to denote partial derivatives) $\displaystyle{\det\left(\begin{array}{cc}u_s & u_t \\ v_s
& v_t\end{array}\right)=u_sv_t-u_tv_s=0}$. Standard recognition of singularities gives:\\
(a) $\Sigma$ is empty\\
(b) $\Sigma$ is a smooth curve and the restriction $g|\Sigma$ of $g$ to $\Sigma$ is nonsingular\\
(c) $\Sigma$ is a smooth curve and $g|\Sigma$ has a cusp of the form $t\to (t^2,t^3)$\\
(d) $\Sigma$ is an isolated point for a lips and a pair of transverse curves for a beaks\\
(e) $\Sigma$ is a smooth curve and $g|\Sigma$ is of the form $t\to (t^3, t^4)$.

\medskip

Ideally we would like a direct link between the contact of a sphere and the nature of the mapping
$g$, but as yet we do not know of such a link. Our result is proved on a case by case basis
by direct calculation. We give a brief indication of the calculations in \S\ref{ss:calcs3D}.
We shall see in \S\ref{ss:ondiag3D} that much the same result holds in the on-diagonal case, though
the mapping $g$ is replaced by a different one.

\vspace*{-0.1in}

\subsection{Some calculations}\label{ss:calcs3D}
Starting from the equation $G=0$ where $G$ is given in (\ref{eq:pre-SS}), we can differentiate with respect
to $s$ and $t$ using standard results about the derivatives of tangent vectors and normal vectors along
curves on a surface; see for example~\cite[\S6.1]{solidshape} for the relevant formulas. Needless to say the
details become complicated by the time the second or third derivatives are reached, and we suppress these
details here, merely giving an indication of how the calculations begin.

Differentiating $G=0$ with respect to $s$ and $t$, which are arclengths on the principal curve $t=0$, $s=0$ respectively,
gives
\begin{eqnarray}
\textstyle{\frac{\partial }{\partial s}}&:& \ee_1(1-r\kappa_1)+r_s(\mm-\nn)=\ff_1(1-r\lambda_1)u_s+\ff_2(1-r\lambda_2)v_s,
\label{eq:Gx}\\
\textstyle{\frac{\partial }{\partial t}}&:& \ee_1(1-r\kappa_2)+r_t(\mm-\nn)=\ff_1(1-r\lambda_1)u_t+\ff_2(1-r\lambda_2)v_t,
\label{eq:Gy}
\end{eqnarray}
valid along the curves $t=0$, $s=0$ respectively.  At points away from these curve we need to allow for the non-unit
speed nature of the parameter curves; however these speeds do not affect our calculations when we evaluate at $s=t=0$.

Consider the $A_1A_1$ case. Then $1-r\kappa_1, \ 1-r\kappa_2, \ 1-r\lambda_1, \ 1-r\lambda_2$ are all nonzero at
$s=t=0$. Taking the vector product of the right-hand sides of (\ref{eq:Gx}) and (\ref{eq:Gy}) then gives 
$(1-r\lambda_1)(1-r\lambda_2)(u_sv_t-v_su_t)\nn$. The mapping $g: (s,t)\to(u,v)$ is a local diffeomorphism at $s=t=0$ if and only if this is nonzero.
However, if the left-hand sides of   (\ref{eq:Gx}) and (\ref{eq:Gy}) are parallel vectors, then $\ee_1, \ee_2$ and $\mm-\nn$
are linearly dependent, and this is impossible since $\ee_1, \ee_2$ are tangent vectors to the sphere $S_0$,
and $\mm-\nn$ is parallel to the chord $\pp-\qq$ of the same sphere. Hence in the $A_1A_1$ case, $g$ is a local
diffeomorphism, verifying assertion (a) of Proposition~\ref{prop:AkA1cases}.

Keeping the contact at $\qq_0$ as $A_1$, so that $1-r\lambda_1, \ 1-r\lambda_2$ are nonzero, but moving to
$A_2$ contact at $\pp_0$ we have $1-r\kappa_1=0$; note that for $A_3$ (ridge or crest point) we will have also $\kappa_{1s}=0$:
the derivative of $\kappa_1$ in the first principal direction equal to zero~\cite[p.144]{mumfordbook}.
Using the linear independence of $\mm-\nn, \ \ff_1, \ff_2$ (\ref{eq:Gx}) now gives $r_s, u_s, v_s$ all equal to 0 at $s=t=0$,
while $r_y, u_y, v_y$ can all be evaluated at $s=t=0$ from (\ref{eq:Gy}). For example,
$r_y=-\ee_2\cdot\nn(1-r\kappa_2)/(\mm\cdot\nn-1)$, where the denominator cannot be 0 since the points
$\pp_0, \qq_0$ are distinct. We can now check the assertion (b) of Proposition~\ref{prop:AkA1cases}.

Continuing in this way we find, for $A_3A_1$, that $r_{ss}, u_{ss}, t_{ss}$ are all zero at $s=t=0$, and the
condition for $\Sigma$ to be singular comes to $u_{st}v_t-u_tv_{st}=0$. A very similar argument to that in the
$A_1A_1$ case enables us to express this condition as $\ee_2\cdot\nn(1-r\kappa_2)\kappa_1^2=\kappa_{1t}(\mm\cdot\nn-1)$,
and after some manipulation this corresponds precisely to the transitional $A_3A_1$ condition found by
Pollitt~\cite{pollitt-thesis}: more picturesquely, it means that the osculating plane of the line of curvature
on $M$ at $\pp_0$ in the direction $\ee_1$ passes through $\qq_0$. In this way we verify (c) and (d)
of Proposition~\ref{prop:AkA1cases}. 

Unfortunately the criterion which distinguishes the two $A_3A_1$
transitional cases, and which therefore distinguishes lips from beaks, though computable, is complicated
and we do not have a simple geometrical interpretation of it. There is a strong link with the function
$R(s,t)=r\kappa_1$: in fact in the transitional $A_3A_1$ situation, $R_s=R_t=0$ at $s=t=0$, and
the `lips' case corresponds to $R$ having a maximum or minimum and the `beaks' case to $R$ having a saddle point.

\vspace*{-0.1in}
\subsection{The local (on-diagonal) case}\label{ss:ondiag3D}

Here we consider the case of the pre-symmetry set of a single surface piece $M$ corresponding to points close to an $A_{\ge 3}$ singularity,
at $\pp_0$ say, which may be taken as the origin in $\RR^3$. This case has
the additional difficulty that the points in question intersect the diagonal of the space $\RR^4$, that is the
set of points $(s,t,s,t)$. This means that, in the notation of (\ref{eq:pre-SS}), $G^{-1}({\bf 0})$ will never be smooth
and we have to `eliminate' the diagonal component in some way. It turns out that we can often do this
by using a different parametrization from that in \S\ref{ss:offdiag3D}. We give brief details below.

As before, we let $(s,t)$ and $(u,v)$ be the parameters for points of contact of a bitangent sphere,
where now all four numbers $s,t,u,v$ are close to 0. In practice we can take $M$ in Monge
form $z=f(x,y)$ so that the contact points are $\pp=(s,t,f(s,t))$ and $\qq=(u,v,f(u,v))$:
\begin{eqnarray}
f(x,y)&=&\textstyle{\frac{1}{2}}\displaystyle (\kappa_1x^2+\kappa_2y^2)+b_0x^3+b_1x^2y+b_2xy^2+b_3y^3 \nonumber \\
&& + c_0x^4 + \ldots + c_4y^4 + d_0x^5 + \ldots +d_5y^5 + \ldots.
\label{eq:monge}
\end{eqnarray}
Instead of having $r,u,v$ expressed as smooth functions of $s,t$ it is much better to use
$s,u$ as the parameters and express $r,t,v$ in terms of these. The great advantage of this is that
the diagonal points appear as $s=u$, that is, the intersection of the diagonal with the true
pre-symmetry set appears as a {\em curve} $\{ (s,t(s,s),s,v(s,s)\}$ on the pre-symmetry set (pre-SS).

The pre-SS is symmetrical about the diagonal in the sense that, for all $a,b,c,d$, 
$(a,b,c,d)\in \mbox{pre-SS} \Longleftrightarrow (c,d,a,b) \in \mbox{pre-SS}$.
This implies that, 
\begin{equation}
\mbox{for all} \  (s,u) \in \RR^2 \ \mbox{we have} \  t(s,u) = v(u,s).
\label{eq:symm}
\end{equation}

We find the following.  Let $h:(s,u)\to (s,t(s,u))$ be the mapping determined by $t$ as a function 
of $s$ and $u$. The relation between this and the mapping $(s,u)\to (u,v(s,u))$ is completely
symmetric so we need consider only $h$. The role of $h$ is to parametrize the points \pp\ 
close to the origin which are one point of contact of a bitangent sphere. The other point
of contact is parametrized by $(u,v(s,u))$.
\begin{proposition}
\begin{enumerate}
\item[{\rm (a)}] $A_3$: $h$ is a local diffeomorphism (see note below),
\item[{\rm (b)}] $A_4$: $h$ is a fold mapping,
\item[{\rm (c)}] $A_5$: $h$ is a lips mapping (beaks mappings do not occur).
\end{enumerate}
\label{prop:ondiagonal}
\end{proposition}
In all cases the pre-symmetry set is the graph of a smooth function $h$ and is therefore itself smooth. Fig.~\ref{fig:projs} illustrates
these cases (see also Fig.~\ref{fig:fold}). \\
{\bf Note} on $A_3$. There is an exception, which occurs when the tangent to the ridge (crest line) through
the origin is in the `other' principal direction. That is, if the ridge corresponds to the principal
direction $\ee_1$ then the tangent is in the direction $\ee_2$. This is slightly mysterious, but in this
case we can use a less satisfactory parametrization, by $s$ and $t$, instead.

There follow some brief details of the above cases.

\smallskip\noindent
${\mathbf A_3}$: We find
\begin{equation}
t(s,u)=\alpha(s+u) + \ \mbox{h.o.t.}, \ \ v(s,u) = \alpha(s+u) + \ \mbox{h.o.t.},
\label{eq:A3local}
\end{equation}
and $\alpha = \frac{1}{4}(\kappa_1^3\kappa_2-\kappa_1^4-8c_0\kappa_2+8c_0\kappa_1+4b_1^2)/(c_1\kappa_2-c_1\kappa_1-2b_1b_2)$.
The numerator of $\alpha$ is zero precisely when the contact is $A_4$; see below for this case. The denominator being
zero is the exceptional case noted above.

\smallskip\noindent
${\mathbf A_4}$: We find
\begin{equation}
 t(s,u) =t_{20}s^2 + t_{11}su + t_{02}u^2 + \ldots; v(s,u) = t_{02}s^2 + t_{11}su + t_{20}u^2 + \ldots,
\label{eq:A4t}
\end{equation}
and the constant coefficients $t_{ij}$ are determined by the local geometry of the surface, namely
\begin{eqnarray*}
t_{02}&=&\frac{3(d_0\kappa_1^2-2d_0\kappa_1\kappa_2+b_1c_1\kappa_1-b_1c_1\kappa_2+b_1^2b_2+d_0\kappa_2^2)}
{(\kappa_1-\kappa_2)(c_1\kappa_2-c_1\kappa_1-2b_1b_2)}\\
t_{20}&=&t_{02}+b_1/(\kappa_1-\kappa_2), \ \ t_{11}=\textstyle{\frac{4}{3}}t_{02}.
\end{eqnarray*}
The coefficients in $v$ as compared with $t$ are explained by (\ref{eq:symm}). (As with $A_3$ there is an exceptional
case, if the denominator of $t_{02}$ is zero. Here, this means that {\em the ridge is singular}: this
is a transitional $A_4$ case, where the ridge is undergoing a transition. See \cite[p.174]{mumfordbook}.)
The coefficient $t_{02}$ of $u^2$ is zero if and only if the singularity is $A_5$ or higher; thus for $A_4$ it is nonzero. 
From this it is easy to check that the mapping $h:(s,u)\to(s,t)$ is a fold. The critical set $\Sigma$ of $h$
is given by $t_{11}s +2t_{02}u + \ldots  = 0$ and has tangent line $2s+3u=0$.
The fold mapping ensures that, given a point  $(s_1,u_1)$ close to the origin, and not on $\Sigma$, there is
a second distinct point $(s_2,u_2)$ with the same image under $\gamma$. This means that $s_1=s_2=s$ say, and
$t(s,u_1)=t(s,u_2)=t$, say.  Then $(s,t,u_1,v(s,u_1))$ and $(s,t,u_2,v(s,u_2))$ both belong to the pre-symmetry set,
and the points $(u_1,v(s,u_1)), (u_2,v(s,u_2))$ are the distinct points of contact of two 
spheres, both tangent
to $M$ at $(s,t)$. Compare Figure~\ref{fig:fold}.

\smallskip\noindent
${\mathbf A_5}$ We find
\begin{equation}
t = \frac{b_1}{\kappa_1-\kappa_2}s^2+ t_{30}s^3+t_{21}s^2u+t_{21}su^2+\textstyle{\frac{2}{3}}\displaystyle t_{21}u^3 + \mbox{h.o.t}.
\label{eq:A5t}
\end{equation} 
The critical set of the mapping $(s,u)\mapsto (s,t)$ is now $t_{21}s^2+2t_{21}su+2t_{21}u^2+$ h.o.t. = 0. The discriminant
of this quadratic form is $-4t_{21}^2<0$ so that the critical set of the mapping has an isolated point at the origin.
This means that it can be only `lips' and not `beaks'.

\section{Conclusion}\label{s:conclusion}
\vspace*{-0.1in}
We have studied pre-symmetry sets of surfaces in 3D, analysing many of them as graphs of functions. This approach brings out
the geometry of the pre-symmetry sets and allows one to see them evolving through a generic transition. There are many other interesting
cases to study, particularly those where the pre-symmetry set is in fact a singular surface.

\smallskip\noindent\begin{small}
{\em Acknowledgements:} 
This work is a part of the DSSCV project supported by the IST Programme of the European Union (IST-2001-35443).
The first author was supported by this grant. The authors are also grateful
to V. Zakalyukin, A. Pollitt and A. Kuijper for helpful discussions.
\end{small}

\end{document}